\documentclass[12pt]{amsart}
\usepackage{amssymb}
\usepackage{mathtools}
\usepackage[pdftex]{hyperref}
\input xy
\xyoption{all}

\usepackage{mathrsfs}

\hypersetup{colorlinks,citecolor=blue,linktocpage}
\usepackage{amsmath,amscd,amssymb,euscript}
\usepackage{mathtools}
\usepackage{latexsym}
\usepackage{graphicx}
\usepackage{tikz}
\usetikzlibrary{backgrounds}
\usetikzlibrary{decorations.fractals}
\usetikzlibrary{calc,intersections,through,backgrounds}
\usetikzlibrary{arrows.meta}
\tikzset{In/.tip = {Hooks[right]}}
\tikzset{Onto/.tip = {To[sep] To}}
\tikzset{Eq/.style = {-, double equal sign distance}}
\tikzset{Iso/.style = {edge node = {node[above] {$\sim$}}, inner sep = 0}}
\usetikzlibrary{positioning}
\usepackage{mathrsfs}
\usepackage{epstopdf}
\usepackage{rotating}
\usepackage{lscape}
\usepackage{color}
\definecolor{alizarin}{rgb}{0.82, 0.1, 0.26}
\definecolor{darkred}{rgb}{0.55, 0.0, 0.0}
\definecolor{brightmaroon}{rgb}{0.76, 0.13, 0.28}

\newcommand{\del}{\partial}
\newcommand{\odel}{\overline{\partial}}

\newcommand{\bC}{{\mathbb C}}

\newcommand{\bG}{{\mathbb G}}
\newcommand{\bH}{{\mathbb H}}

\newcommand{\bP}{{\mathbb P}}
\newcommand{\bQ}{{\mathbb Q}}
\newcommand{\bR}{{\mathbb R}}
\newcommand{\bS}{{\mathbb S}}

\newcommand{\bV}{{\mathbb V}}
\newcommand{\bZ}{{\mathbb Z}}
\newcommand{\cA}{{\mathcal A}}

\newcommand{\cD}{{\mathcal D}}
\newcommand{\cE}{{\mathcal E}}
\newcommand{\cF}{{\mathcal F}}
\newcommand{\cG}{{\mathcal G}}
\newcommand{\cH}{{\mathcal H}}

\newcommand{\cJ}{{\mathcal J}}

\newcommand{\cM}{{\mathcal M}}

\newcommand{\cO}{{\mathcal O}}

\newcommand{\cU}{{\mathcal U}}
\newcommand{\cV}{{\mathcal V}}

\newcommand{\cX}{{\mathcal X}}
\newcommand{\cY}{{\mathcal Y}}
\newcommand{\cZ}{{\mathcal Z}}

\newcommand{\oS}{\overline{S}}

\newcommand{\ov}{\overline{v}}

\newcommand{\oz}{\overline{z}}
\newcommand{\oalpha}{\overline{\alpha}}
\newcommand{\onabla}{\overline{\nabla}}
\newcommand{\osigma}{\overline{\sigma}}

\newcommand{\tnu}{\widetilde{\nu}}

\newcommand{\tX}{\widetilde{X}}

\newcommand{\ra}{\rightarrow}
\newcommand{\lra}{\longrightarrow}

\newcommand{\surj}{\mathrel{\mathrlap{\rightarrow}\mkern1mu\rightarrow}}

\newtheorem{proposition}{Proposition}[section]
\newtheorem{lemma}[proposition]{Lemma}

\newtheorem{definition}[proposition]{Definition}

\theoremstyle{definition}

\numberwithin{equation}{section}

\newcommand{\Coker}{\operatorname{Coker}}

\newcommand{\Ext}{\operatorname{Ext}}
\newcommand{\GL}{\operatorname{GL}}
\newcommand{\Grass}{\operatorname{Grass}}

\newcommand{\Hom}{\operatorname{Hom}}
\newcommand{\Id}{\operatorname{Id}}
\newcommand{\Image}{\operatorname{Im}}
\newcommand{\Ker}{\operatorname{Ker}}

\begin{document}

\title{On infinitesimal invariants of normal functions}

\author{Elham Izadi}

\address{Department of Mathematics, University of California San Diego, 9500 Gilman Drive \# 0112, La Jolla, CA
92093-0112, USA}

\email{eizadi@math.ucsd.edu}

\dedicatory{To the memory of Alberto Collino}

\thanks{}


\begin{abstract}
We present an overview of some of Alberto Collino's work which uses the Griffiths infinitesimal invariant of a normal function.
\end{abstract}

\maketitle

\tableofcontents

\section*{Introduction}

Let $X$ be a smooth complex projective variety. For divisors on $X$, it is well-known that algebraic and homological equivalence coincide. Until the work of Griffiths however, it was not clear whether the same would hold for cycles of higher codimension. Griffiths found the first examples of algebraic varieties for which homological and algebraic equivalence were different \cite{Griffiths1969}. Accordingly, the group of homologically trivial cycles modulo algebraic equivalence is called the Griffiths group of $X$. 

To study algebraic and homological equivalence and to generalize Lefschetz' proof of the $(1,1)$-theorem to higher codimension cycles, Griffiths introduced the Griffiths intermediate jacobians $J^d(X)$, correcting an earlier notion introduced by Weil. A homologically trivial algebraic cycle $Z$ of codimension $d$ has an Abel-Jacobi image $AJ(Z)\in J^d (X)$ obtained, roughly speaking, via integration of differential forms on a topological chain with boundary $Z$. The intermediate jacobian is a complex torus and it has an algebraic part $J^d(X)_a$ which is the largest abelian variety contained in it. It is not too difficult to see \cite[Section 13]{Griffiths1969}, that if $Z$ is algebraically equivalent to $0$, then $AJ(Z) \in J(X)_a$.

Families of cycles (and primitive cohomology classes) give rise to sections of families of intermediate jacobians called normal functions. The infinitesimal invariant of a normal function is, roughly speaking, the family of its derivatives. The Griffiths infinitesimal invariant was refined by Green \cite{Green1989-Inf} and used to obtain many interesting geometric results by Collino, Green, Pirola, Voisin.

One such set of results has to do with the cycle $C - C^-$ (defined up to translation) in the Jacobian of a smooth curve $C$. Ceresa \cite{Ceresa1983} proved earlier that this cycle is not algebraically equivalent to $0$ for a very general curve. Using the infinitesimal invariant of the normal function of $C - C^-$, Collino and Pirola \cite{CollinoPirola1995} gave a second proof of Ceresa's result. They also showed that, in genus $3$, this infinitesimal invariant could be identified with the equation of the canonical curve when $C$ is not hyperelliptic. They further showed that Ceresa's theorem also holds for very general plane curves of degree $\geq 4$ and also for very general curves in subvarieties of codimension $< \frac{g+2}{3}$ of $\cM_g$ if $g>3$. In genus $3$, they showed that the Ceresa cycle is not algebraically equivalent to $0$ on any subvariety of codimension at most $2$ of $\cM_3$ which is not contained in the hyperelliptic locus.

Collino and Pirola \cite{CollinoPirola1995} introduced an adjunction map for infinitesimal deformations of smooth varieties. They used this map, in conjunction with Griffiths' formula in Section \ref{secGriffform}, to compute certain values of the infinitesimal invariant of the normal function of $C - C^-$. They showed that the normal function is non-abelian, meaning it does not map into the abelian part of the intermediate jacobian. This implies that the Ceresa cycle is not algebraically equivalent to $0$ for a very general curve.

Pirola and Zucconi generalized the Collino-Pirola constructions to produce a bound for the geometric genus of subvarieties of general type of generic abelian varieties. For instance, they proved that if $X$ has dimension $n$, geometric genus $p_g (X)$, and admits a generically finite morphism to a generic abelian variety of dimension $a$, then
\[
p_q (X) \geq {a-n \choose 2} + {a \choose n}.
\]
When $X$ is a surface of general type whose Albanese map is generically injective, they proved Castelnuovo's inequality:
\[
m \leq p_g (X) +2q -3,
\]
where $q$ is the irregularity of $X$ (i.e., the dimension of its Albanese variety) and $m$ is the maximum variation of $X$, meaning the dimension of a family of deformations of $X$ mapping to abelian varieties of dimension $a$ whose Kodaira--Spencer map is injective.
They also generalized the results of Collino--Pirola by showing that the only subvariety of dimension at least $2g-1$ of $\cM_g$ on which the Abel-Jacobi image of the Ceresa cycle is trivial is the hyperelliptic locus.

Smooth cubic threefolds provided one of the first examples of non-rational unirational threefolds \cite{ClemensGriffiths1972}. The intermediate jacobian $J(X)$ of a cubic threefold $X$ is isomorphic to the Albanese variety of its Fano variety of lines $F$ via the Abel-Jacobi map. Via this isomorphism, the Fano variety $F$ admits embeddings into $J(X)$ and has, similarly to the curve case, a cycle $F - F^-$ well-defined up to translation. Collino, Naranjo and Pirola \cite{CollinoNaranjoPirola} computed the infinitesimal invariant of the normal function of this cycle. They used it to show that, for $X$ general, the cycle $F - F^-$ is not algebraically equivalent to $0$ and that there is no divisor in $J(X)$ containing both $F$ and $F^-$ where $F$ and $F^-$ are homologically equivalent. They further showed that the infinitesimal invariant determines the curve of double lines in $F$ which in turn determines the cubic threefold $X$.

We give an overview of the tools and methods used to obtain the above results. We begin by giving the basic definitions of Hodge structures (Section \ref{defHodge}), period domains (Section \ref{defperiod}) and variations of Hodge structures (Sections \ref{defhor} and \ref{defVHS}). Next we define the Griffiths intermediate jacobian and give an overview of some important results about the Griffiths group (Section \ref{sectIntJ}). We define normal functions (Section \ref{defnorm}), infinitesimal variations of Hodge structure (Section \ref{secInfVar}) and the infinitesimal invariant of a normal function (Section \ref{SectInf}). We explain Griffiths' formula for the calculation of the infinitesimal invariant in Section \ref{secGriffform} and Voisin's formula in Section \ref{sectvoisin}. In Section \ref{defadj} we define the adjunction map and in Section \ref{sectlink} we explain how the adjunction map was used to compute the infinitesimal invariant.

\section{Hodge structures}\label{defHodge}

All schemes are over the field of complex numbers. A variety is a reduced separated scheme of finite type over $\bC$. In what follows, $X$ will denote a smooth projective variety of dimension $n$. We begin with some definitions and explaining our setup.

\begin{definition}
A rational (respectively, integral) Hodge structure of weight $k \in \bZ$ is the datum of a
$\bQ$-vector space (respectively, free $\bZ$-module) $V$ such that the complexification
\[
V_{\bC} := V\otimes \bC
\]
admits a direct sum decomposition
\[
V_{\bC}  = \oplus_{p+q = k} V^{p,q} \quad \hbox{   with   } \quad \overline{V^{p,q}} = V^{q,p} \: (p,q \in \bZ),
\]
where the bar indicates complex conjugation defined by $\overline{v\otimes z} = v \otimes
\overline{z}$ for $v\in V, z\in \bC$. A Hodge structure (of positive weight) is called
effective if $V^{p,q} =0$ when either $p$ or $q$ is negative.
\end{definition}

\begin{definition}

A polarization on the rational (respectively, integral) Hodge structure $V$ is the datum of a bilinear map
\[
\psi : V\times V \lra \bQ \quad (\hbox{respectively, }V\times V \lra \bZ),
\]
symmetric for even $k$, skew-symmetric for odd $k$, whose complexification (i.e., linear
extension to $V_{\bC}$) satisfies
\begin{equation}\label{Riem1}
\psi (H^{p,q}, H^{r,s} ) = 0 \quad \hbox{ unless } \quad p=s, q=r
\end{equation}
and
\begin{equation}\label{Riem2}
i^{p-q} \psi (v, \ov ) > 0 \quad \hbox{ for } \quad v \in H^{p,q}, v\neq 0.
\end{equation}
In particular, $\psi$ is always nondegenerate.\\
Equations \eqref{Riem1} and \eqref{Riem2} are referred to, respectively, as the first and the second Hodge-Riemann bilinear relations.
\end{definition}

Given a smooth complex projective variety $X$, the rational cohomology $V := H^k (X, \bQ)$ (respectively, integral cohomology modulo torsion)
carries a natural Hodge structure:
\[
V\otimes \bC = H^k (X, \bQ )\otimes \bC \cong H^k (X, \bC) \cong \oplus H^{p,q} (X) =
\oplus V^{p,q}
\]
where $V^{p,q} = H^{p,q} (X)$ is, as usual, the space of cohomology classes of differential forms of
type $(p,q)$. To obtain polarizations on the cohomology of algebraic varieties, we need to restrict ourselves to the primitive cohomology.

Given an ample line bundle $L$ on $X$, let $\eta := c_1 (L) \in H^2 (X,
\bZ)$ be its topological first Chern class. Left multiplication by $\eta$ defines a linear
map
\[
\Lambda : H^k (X, \bZ ) \lra H^{k+2 } (V, \bZ)
\]
called the Lefschetz operator.
For $0\leq k\leq n:= \hbox{dim}_{\bC} X$, the $(n-k)$-th power of $\Lambda$ is
injective. The primitive part $P^k (X, \bZ)$ of $H^k (X, \bZ)$ (respectively, $P^k (X, \bC)$
of $H^k (X, \bC)$) is defined to be the kernel of
$\Lambda^{n-k+1}$. If $k > n$, the primitive part is defined to be $0$. The primitive
cohomology carries a structure of polarized rational (or integral) Hodge structure. The spaces $P^{p,q}$
are defined to be the intersections
\[
P^{p,q} := P^k \cap H^{p,q}
\]
and the polarization is the Hodge bilinear form defined by
\[
\psi (v,w ) = (-1)^{\frac{k(k-1)}{2}} (\Lambda^{n-k} v \wedge w) [X] = (-1)^{\frac{k(k-1)}{2}} \int_X \eta^{n-k} \wedge v \wedge w,
\]
where $[X]$ is the fundamental class of $X$. The Hodge-Riemann bilinear relations assert that this
is indeed a polarization in the sense that we defined earlier (see, e.g.,  \cite[Vol. I, Chapter 6]{Voisin2003-Hodge}).

The Hard Lefschetz Theorem is the assertion that one has the decomposition
\[
H^k (X,\bQ ) = \oplus \Lambda^j P^{k-2j} (X, \bQ), \quad j \geq max\{k-n, 0\}.
\]
So the primitive cohomology determines the cohomology of an algebraic variety.

To construct a parameter space for Hodge structures, one cannot use the above definition
since the $H^{p,q}$ do not vary holomorphically on a family of algebraic varieties. Fortunately, there
is an equivalent definition of a Hodge structure that allows one to construct
a complex analytic parameter space, called a period domain.

\begin{definition}
A rational (respectively, integral) Hodge structure of weight $k$ is the data of a $\bQ$ vector space $V$ (respectively, free $\bZ$-module) together
with a decreasing filtration $F^p$ on its complexification $V_{\bC}$ such that
\[
\forall p \quad V_{\bC} = F^p \oplus \overline{F^{k-p+1}}.
\]

\end{definition}

The relation between the two definitions of Hodge structure is
\[
F^p := \oplus_{i\geq p} V^{i, k-i}
\]
and
\[
V^{p,q} = F^p \cap \overline{F^q}.
\]
In terms of the Hodge filtration, the Hodge--Riemann bilinear relations become
\[
\psi (F^p, F^{k-p+1}) = 0 \quad \hbox{ for all } p
\]
and
\[
\psi (Cv, \ov) > 0 \hbox{ for } v\in V_{\bC}, v\neq 0
\]
where $C$ is the Weil operator defined by $Cv = i^{p-q} v$ if $v\in V^{p,q}$.

Griffiths proved that the Hodge filtration varies holomorphically in
families (algebraically for a family of algebraic varieties, see \cite{Griffiths1968-II}).

Yet a third definition of a Hodge structure was given by Deligne \cite{Deligne1972-II} using representations. This point of view allowed the introduction of Mumford--Tate groups which have been very useful in the study of Hodge structures and period domains, see, e.g., \cite{GreenGriffithsKerr2012}.

Denote by $\bS$ the group scheme whose set of points $\bS (A)$ on an algebra $A$ is the subgroup of $\GL_2 (A)$ of matrices
\[
\left(\begin{array}{cc}
a & b \\
-b & a
\end{array}\right).
\]
The group $\bS$ naturally contains $\bG_m$ as the subgroup of diagonal matrices and maps to $\bG_m$ via the norm map
\[
\left(\begin{array}{cc}
a & b \\
-b & a
\end{array}\right) \longmapsto a^2 + b^2.
\]
The group $\bS (\bR)$ is naturally isomorphic to $\bC^*$ via the map
\[
a+ib \longmapsto  \left(\begin{array}{cc}
a & b \\
-b & a
\end{array}\right).
\]
\begin{definition}
A rational (respectively, integral) Hodge structure of weight $k$ is the data of a $\bQ$ vector space (respectively, free $\bZ$-module) $V$ together
with a representation
\[
\bS (\bR) \lra \GL (V\otimes \bR).
\]
\end{definition}
Under the identification $\bS (\bR) = \bC^*$, the Hodge component $V^{p,q}$ is the summand of $V\otimes \bC$ where $z$ acts as multiplication by $z^p \oz^q$.

\section{Period domains}\label{defperiod}

For the results of this section, we refer to \cite{Griffiths1968-I}, \cite{Griffiths1968-II}, \cite{GriffithsSchmid1969} and \cite{Schmid1973}. Fix a finitely generated free abelian group $V$ with complexification $V_{\bC}$. Also fix
an integer $k$ and a collection of non-negative integers $\{ h^{p,q}\}$ which satisfy
$h^{p,q} = h^{q,p}$, $\sum h^{p,q} = d := \hbox{rank}_{\bZ}V$ and $h^{p,q} \neq 0$ only when
$p+q =k$.

Let $\check{\cF}$ be the flag variety of filtrations $F^\bullet$ of $V_{\bC}$ of dimensions $f_p := \dim F^p =
\sum_{i\geq p} h^{i, k-i}$. Then $\check{\cF}$ can be realized as a closed subvariety of a
product of Grassmannians which shows that $\check{\cF}$ is a projective algebraic
variety. Furthermore $GL (V_{\bC})$ acts algebraically and transitively on
$\check{\cF}$ which shows, in particular, that $\check{\cF}$ is smooth. The filtrations
that satisfy
\[
V_{\bC} = F^p \oplus \overline{F^{k-p+1}}
\]
form an open subset $\cF$ of $\check{\cF}$ in the complex analytic topology. Therefore
$\cF$ is a complex manifold that parametrizes Hodge structures of weight $k$ on $V$ with
Hodge numbers $h^{p,q}$.

Now let $\psi$ be a nondegenerate bilinear form on $V$ (or $V_{\bQ}$, if one wishes to work over the rationals), symmetric if $k$ is even,
skew-symmetric if $k$ is odd. Let $\check{\cD} \subset \check{\cF}$ be the closed
algebraic subvariety of those filtrations that satisfy the first Hodge-Riemann bilinear relation \eqref{Riem1}:
\[
\psi (F^p, F^{k-p+1}) = 0 \quad \hbox{ for all } p.
\]
The orthogonal group of the bilinear form $\psi$ acts transitively on $\check{\cD}$, hence
$\check{\cD}$ is a smooth projective algebraic variety. In fact, if we denote $G_\bR$ the group of automorphisms of $V_\bR$ that preserve $\psi$, then $\check{\cD} = G_\bC/P$ for a parabolic subgroup $P$ of $G_\bC$ ($G_\bC$ is the base change of $G_\bR$ to $\bC$). The filtrations in $\check{\cD}$
that satisfy the second Riemann bilinear relation \eqref{Riem2}
\[
\psi (Cv, \ov) > 0 \hbox{ for } v\in V_{\bC}, v\neq 0
\]
form an open subset $\cD$ of $\check{\cD}$ in the complex analytic topology, hence a
complex manifold. One can show that the group $K := P\cap G_\bR$ is compact and, under the action of $G_\bR$, $D$ is identifed with the homogeneous complex manifold $G_\bR/K$.

Furthermore, one can show that the subgroup $G_\bZ$ of $G_\bR$ of elements that
leave the lattice $V$ globally invariant acts properly discontinuously on $\cD$. Hence the
quotient $G_\bZ \backslash \cD$ is a complex analytic variety. This quotient is called the period
domain. The analogous statement about classifying spaces for weighted hodge structures
without polarization fails. This is one of the reasons why one considers polarized Hodge
structures.

Given a smooth projective family of algebraic varieties $f : \cX \ra S$, let $s_0\in S$ be a base point and $X := X_{s_0}$ the fiber of $f$ at $s_0$. Then, choosing $V := P^k (X, \bZ)$, the assignment $s\mapsto (P^k (X_s, \bZ), \psi_s)$ defines a holomorphic map $S\mapsto D / \Gamma$. This map is locally liftable to $D$ and its differential maps the space $T_sS$ into a subspace of the tangent space of $D$ called the horizontal tangent space \cite{CornalbaGriffiths1975}, \cite{Griffiths1969} which we define in Section \ref{defhor}.

\section{The horizontal tangent bundle}\label{defhor}

Under the embedding
\[
\check{\cD} \subset \check{\cF} \subset \prod \Grass_p (f_p, d),
\]
we have an embedding of the tangent space to $\check{D}$ at a point $F = (F^p)$ in the tangent space to the product of Grassmannians:
\[
T_F \check{D} \subset \oplus_p \Hom (F^p , V_\bC / F^p ).
\]
One easily checks that $T_F \check{D}$ is the set of $\xi = \oplus_p \xi_p \in \oplus_p \Hom (F^p , V_\bC / F^p )$ satisfying the following conditions.
\begin{enumerate}
\item The diagram
\[
\xymatrix{
F^p \ar[d] \ar[r]^{\xi_p} & V_\bC/ F^p \ar[d] \\
F^{p-1} \ar[r]^{\xi_{p-1}\quad} & \quad V_\bC/F^{p-1}}
\]
is commutative
\item and
\[
\psi (\xi_p v, w ) + \psi (v, \xi_{n-p+1} w ) =0, \hbox{ for all } v\in F^p, w\in F^{n-p+1}.
\]
\end{enumerate}

\begin{definition}
The horizontal tangent space $T_{h,F}$ to $\check{D}$ is the set of vectors $\xi = (\xi_p) \in T_F \check{D}$ satisfying the infinitesimal period relations or Griffiths transversality
\[
\xi (F^p) \subset F^{p-1}
\]
for all $p$.
These tangent spaces form the horizontal tangent bundle $T_{h,\check{D}} \subset T_{\check{D}}$. The horizontal tangent bundle of $D$ is the restriction $T_{h, D} := T_{h,\check{D}} |_D$.
\end{definition}
It is clear that the horizontal tangent bundle is a holomorphic subbundle of $T_{\check{D}}$, invariant under the action of $G_\bC$.

Note that an element $\xi = (\xi_p) \in T_{h,F} \check{D}$ induces maps
\[
V^{p,q} = F^p /F^{p+1} \lra V^{p-1, q+1} = F^{p-1} / F^p
\]
for all $p,q$. We therefore have a natural embedding
\[
T_{h,F} \check{D} \subset \oplus_{1\leq p\leq k} \Hom (V^{p,q} , V^{p-1, q+1}).
\]

\section{Variation of Hodge structure}\label{defVHS}

When properly defined, the datum of a morphism to a period domain is equivalent to that of a variation of Hodge structure. Before giving the definition of a variation of Hodge structure, we recall that of Gauss--Manin connections.

Given a complex manifold $S$ with a local system $\bV$ of finitely generated free abelian groups or rational or real vector spaces, the holomorphic vector bundle $\cV := \bV \otimes \cO_S$ carries a natural flat connection whose local system of flat sections is naturally identified with $\bV$. This is the Gauss--Manin connection 
\[
\nabla : \cV \lra \cV \otimes \Omega^1_S
\]
which can be locally defined as follows.

Choose a local basis $\{\sigma_1, \ldots, \sigma_d\}$ of $\bV$. Then the Gauss--Manin connection sends a local section $\sigma = \sum \alpha_i \sigma_i$ of $\cV$ to $\nabla (\sigma) = \sum \sigma_i \otimes d \alpha_i$. One checks that the definition is independent of the choice of local basis, hence glues to define a global connection and that its curvature form is $0$ (see, e.g., \cite[Vol. I, Chapter 9]{Voisin2003-Hodge}).

\begin{definition}
A variation of Hodge structure $(S, \bV, \cF^\bullet, \nabla)$ of weight $k$ on a complex manifold $S$ is the data of a local system $\bV$ of finitely generated free abelian groups (or finite dimensional rational or real vector spaces) with a decreasing filtration
\[
\{0 \} \subset \cF^k \subset \ldots \subset \cF^0 = \cV
\]
of the holomorphic bundle $\cV := \bV \otimes \cO_S$ such that if $\nabla : \cV \ra \cV \otimes \Omega^1_S$ denotes the Gauss-Manin connection, then we have Griffiths transversality (or the infinitesimal period relation)
\[
\nabla (\cF^p )\subset \cF^{p-1} \otimes \Omega^1_S
\]
and we have $C^\infty$ splittings
\[
\cV = \cF^p \oplus \overline{\cF^{k-p+1}}
\]
for all $p$. The Hodge bundles are then defined as
\[
\cV^{p, q} := \cF^p / \cF^{p+1}.
\]
They satisfy
\[
\cV^{p,q} = \overline {\cV^{p,q}}
\]
and we have a $C^\infty$ splitting
\[
\cV = \oplus_{p+q=k} \cV^{p,q}.
\]
A polarization of an integral or rational variation of Hodge structure is the datum of a {\em flat} nondegenerate bilinear form
\[
\psi : \bV \times \bV \lra \bZ_S \, (\hbox{resp. } \bQ_S)
\]
which satisfies the first and second Hodge-Riemann bilinear relations \eqref{Riem1} and \eqref{Riem2}.
\end{definition}
Note that the flatness of the bilinear form means that it is compatible with the Gauss-Manin connection, i.e., for all sections $s, s'$ of $\cV$,
\[
d(\psi (s, s')) = \psi (\nabla s, s') + \psi (s, \nabla s').
\]
A polarized variation of Hodge structure defines a monodromy representation
\[
\rho : \pi_1 (S, o) \lra G_\bZ
\]
where $o\in S$ is a base point. Let $\Gamma$ be the image of the monodromy representation $\rho$. Sending a point of $S$ to the Hodge structure defined by the polarized variation of Hodge structure defines a period map
\[
S \lra \Gamma \backslash \cD
\]
which is holomorphic, locally liftable and whose local lifts send the tangent bundle $T_S$ into the the horizontal tangent bundle $T_{h,\cD}$. In fact the data of a polarized variation of Hodge structure with monodromy contained in a subgroup $\Gamma$ of $G_\bZ$ is equivalent to the data of a holomorphic map $S \lra \Gamma \backslash \cD$ with the above properties (see \cite{CornalbaGriffiths1975}, \cite{Griffiths1968-II}).

\begin{definition}
An extended variation of Hodge structure is the data of a holomorphic map
\[
f : S \lra \Gamma \backslash \cD
\]
such that the restriction of $f$ to a dense Zariski open subset of $S$ is a variation of Hodge structure.
\end{definition}
Note that the difference between a variation of Hodge structure and an extended variation of Hodge structure is that for an extended variation there exists a proper analytic subset $Z$ of $S$ such that $f$ fails to be locally liftable at the points of $Z$. If $S$ is a Zariski open subset of a smooth projective variety $\oS$, given a polarized variation of Hodge structure $f : S \ra \Gamma \backslash \cD$, there is a maximal Zariski open set $S' \subset \oS$ such that $f$ extends to an extended variation of Hodge structure $f' : S' \ra \Gamma \backslash \cD$ where $f'$ is {\em proper} (see \cite{CGGH}).

Let $\phi : \cM \ra S$ be a holomorphic family of compact K\"ahler manifolds of complex dimension $n$. For any integer $k\in \{0, \ldots, 2n\}$, we have the variation of Hodge structures $(S, \bH^k := R^k\phi_* \bZ_\cM, \cF^\bullet, \nabla)$, where $\bZ_\cM$ is the constant sheaf with group $\bZ$ on $\cM$, $\cF^\bullet$ is the Hodge filtration on $\cH^k := R^k\phi_* \bZ_\cM \otimes_\bZ \cO_S$ which restricts to the Hodge filtration on $H^k (M_s, \bC)$ for all $s$, and $\nabla : \cH^k \ra \cH^k \otimes \Omega^1_S$ is the Gauss--Manin connection obtained by differentiating sections as above.

For a family of projective algebraic varieties $\cX \ra S$ with a relatively ample line bundle, the primitive cohomology groups of the fibers form local systems $\bV^k$ with polarizations
\[
\Psi : \bV^k \times \bV^k \lra \bZ_S
\]
that give polarized variations of Hodge structures $(S, \bV^k, \cF^\bullet, \nabla, \Psi)$.

\section{Griffiths intermediate jacobians}\label{sectIntJ}

Given an integral Hodge structure of odd weight $k= 2l-1$, the direct sum decomposition
\[
V_\bC = F^l V \oplus \overline{F^l V}
\]
implies that
\[
F^l V \cap V_\bR =0.
\]
Hence the projection
\[
V_\bR \lra V_\bC / F^l V
\]
is an isomorphism of real vector spaces and the lattice $V \subset V_\bR$ maps isomorphically to a lattice in $V_\bC / F^l V$.
\begin{definition}
The intermediate jacobian $J(V)$ is the complex torus
\[
J(V) := \frac{V_\bC}{F^l V \oplus V}.
\]
The algebraic part $J(V)_a$ of the intermediate jacobian is the largest complex subtorus whose tangent space is contained in $V^{l-1, l}$. It is the image, in $J(V)$, of the largest Hodge substructure of $V$ contained in $V^{l-1,l} \oplus V^{l, l-1}$.
\end{definition}
Given a compact K\"ahler manifold $M$ of dimension $m$, Poincar\'e duality gives a commutative diagram
\[
\xymatrix{
H^k (M, \bZ) /_{torsion} \ar[d] \ar[r]^= & H_{2m-k} (M, \bZ) /_{torsion} \ar[d] \\
H^k (M, \bC) \ar[r]^= & \quad H^{2m-k}(M, \bC)^*}
\]
where the second vertical map is induced by integration of differential forms on topological cycles. This gives a natural isomorphism
\[
J(V) \cong \frac{F^{m-l+1} H^{2m-k} (M, \bC)^*}{H_{2m-k} (M, \bZ) /_{torsion}}.
\]
Given an algebraic cycle $Z$ of codimension $l$, homologous to $0$, let $W$ be a topological chain with boundary $Z$. One can show (see, e.g., \cite[Vol. I, Chapter 12]{Voisin2003-Hodge}) that integration on $W$ gives a well-defined element
\[
AJ (Z) \in J(V)
\]
called the Abel-Jacobi image of $Z$.

Given a complex manifold $S$, a cycle $\cZ \subset M\times S$ of codimension $l$, flat over $S$ and a point $s_0\in S$, the map
\[
\begin{array}{ccc}
S & \lra & J(V) \\
s & \longmapsto & AJ(Z_s - Z_{s_0}),
\end{array}
\]
where $Z_s$ is the fiber of $\cZ$ at $s$, is holomorphic \cite{Griffiths1968-II}.

It is not too difficult to see, see \cite[Section 13]{Griffiths1969}, that if a cycle $Z$ is algebraically equivalent to $0$, then $AJ(Z) \in J(V)_a$. The group of homologically trivial cycles of codimension $d$ on $M$ modulo algebraic equivalence is called the $d$-th Griffiths group of $M$.

Griffiths proved \cite[Section 14]{Griffiths1969} that the Griffiths groups of generic hyperplane sections of smooth hypersurfaces contain infinite cyclic groups, provided that the original hypersurface carries an algebraic cycle whose homology class is primitive and non-zero. Collino and Ceresa showed \cite{CeresaCollino1983} that in a generic quintic threefold $X$ there are no algebraic equivalence relations between the lines in $X$. In particular, the Abel-Jacobi images of the differences of two distinct lines are non-torsion (by work of Katz \cite{Katz1986}, $M$ contains finitely many lines). Since $J(V)_a=0$ \cite[Section 13]{Griffiths1969}, it follows from the Ceresa--Collino result that the Griffiths group of $X$ is non torsion. Clemens \cite{Clemens1983-GG} proved that the image of the Abel-Jacobi map of a generic quintic threefold is countably generated, and, when tensored with $\bQ$, it is not finitely generated. This implies, in particular, that the Griffiths group tensored with $\bQ$ is not finitely generated. To prove this result, Clemens produced an infinite set of independent elements in the Griffiths group by showing the existence of an infinite sequence of rigid rational curves. Voisin \cite{Voisin1992-GG} gave a second proof of Clemens' theorem where she used the Noether-Lefschetz theorem and the Griffiths infinitesimal invariant of normal functions (see Sections \ref{defnorm} and \ref{SectInf}).

The Lefschetz hyperplane theorem (together with the exponential sequence) implies that the :icard group of $X$ is isomorphic to $\bZ$ if $X$ is a hypersurface in $\bP^{n+1}$ with $n\geq 3$. Griffiths and Harris \cite{GriffithsHarris1985} asked whether the Chow group $CH^2(X)$ of codimension $2$ cycles modulo rational equivalence is isomorphic to $\bZ$ when $X$ is a threefold of degree $d \geq 6$ in $\bP^4$. If this is true then it would follow that the image of the Abel-Jacobi map of $X$ is $0$. Modulo torsion, the latter is a special case of a theorem of Green \cite{Green1989-Inf} and Voisin (unpublished, see \cite[Vol II, Chapter 7]{Voisin2003-Hodge}): they showed that if $X\subset \bP^{2m}$ is a hypersurface of degree $d \geq 2 + \frac{4}{m-1}$, then the image of the $m$-th Abel-Jacobi map of $X$ is torsion.

The main tool in the proof of the Green--Voisin theorem is the Griffiths infinitesimal invariant of a normal function (see Sections \ref{defnorm} and \ref{SectInf}).

The exceptions to the Green--Voisin theorem are the threefolds $X\subset \bP^4$ of degree $d\leq 5$ and cubic fivefolds and sevenfolds. For cubic threefolds, the Abel-Jacobi map from the Albanese variety of the Fano variety of lines in $X$ to the intermediate jacobian of $X$ is an isomorphism \cite{ClemensGriffiths1972}. For quartic threefolds, the Abel-Jacobi map from the Albanese of the family of conics in $X$ to the intermediate jacobian of $X$ is an isomorphism \cite{Letizia1984}. In the case of the cubic fivefold, Collino \cite{Collino1986} showed that the Abel-Jacobi map is an isomorphism between the intermediate jacobian and the Albanese variety of the family of planes in $X$.

In the case of the cubic sevenfold, Albano and Collino \cite{AlbanoCollino1994}, using the techniques we explain below, showed that the Griffiths group of $3$-cycles tensored with $\bQ$ is not finitely generated:

Nori \cite{Nori1993-Conn} proved a far reaching generalization of the Green-Voisin theorem: the Nori connectivity theorem. Nori's theorem implies the Noether--Lefschetz and the Green--Voisin theorems. Nori also defined a natural filtration on the Griffiths group of a smooth projective variety and showed that every graded piece of the filtration can be nontorsion. Using Nori's work \cite{Nori1993-Conn} and Voisin's proof  \cite{Voisin1992-GG} of Clemens' theorem \cite{Clemens1983-GG}, Albano and Collino \cite{AlbanoCollino1994} proved that the Griffiths group of a general cubic sevenfold $M$ is not finitely generated after tensoring with $\bQ$. They also prove the analogous result for the Griffiths group of the complete intersection $Y$ of $M$ with two (or one) generic hypersurfaces of sufficiently high degree. This last result is particularly interesting since, unlike the previous examples, the Griffiths intermediate jacobian of $Y$ is trivial. The nontriviality of the Griffiths intermediate jacobian was previously used in the first proofs of the fact that the Griffiths group of a generic quintic threefold is nontorsion and non finitely generated. Later, Fakhruddin \cite{Fakhruddin1996}, using Nori's connectivity theorem \cite{Nori1993-Conn}, proved that the Griffiths groups of codimension $3$ and $4$ on generic abelian fivefolds are of infinite rank although they are in the kernel of the Abel-Jacobi map.

\section{Normal functions}\label{defnorm}

Given an integral variation of Hodge structure $(S, \bV, \cF^\bullet, \nabla)$ of odd weight $k= 2l-1$ on the complex manifold $S$, one defines the associated family of intermediate jacobians
\[
\cJ_S := \frac{\cV}{\cF^l \oplus \bV}.
\]
Since the Gauss-Manin connection
\[
\nabla : \cV \lra \cV \otimes \Omega^1_S
\]
satisfies
\[
\nabla (\bV) =0, \quad \nabla (\cF^m) \subset \cF^{m-1} \otimes \Omega^1_S,
\]
it induces a map
\[
\onabla : \cJ_S \lra \frac{\cV}{\cF^{l-1}} \otimes \Omega^1_S.
\]
Here we slightly abuse notation by using the same notation for the family of intermediate jacobians and its sheaf of sections.
\begin{definition}
The sheaf of normal functions is
\[
\cJ_{h,S} := \Ker \onabla,
\]
and a normal function is a global section of $\cJ_{h,s}$.
\end{definition}
When $S$ is a quasi-projective variety, a normal function must satisfy additional conditions which govern its behavior at infinity, saying, roughly speaking, that it has at most logarithmic growth (see \cite{ElZeinZucker}).

\section{Infinitesimal variation of Hodge structure}\label{secInfVar}

An infinitesimal variation of Hodge structure is an abstractification of the differential of a variation of Hodge structure. Unlike infinitesimal deformations of algebraic varieties, it contains a little bit more than just first order information.

\begin{definition}
An infinitesimal variation of Hodge structure is the data of a polarized Hodge structure $(V_\bZ, V^{p,q}, \psi)$ together with a vector space $T$ and a linear map
\[
\delta : T \lra \oplus_{1\leq p\leq k} \Hom (V^{p,q}, V^{p-1, p+1})
\]
satisfying the conditions
\begin{equation}\label{infflat}
\delta (\xi_1) \circ \delta (\xi_2) = \delta (\xi_2) \circ \delta (\xi_1),
\end{equation}
\begin{equation}\label{infcomp}
\psi (\delta (\xi) v, w) + \psi(v, \delta (\xi) w) = 0.
\end{equation}
\end{definition}
Condition \eqref{infflat} (which is of second order) expresses the flatness of the Gauss-Manin connection, while condition \eqref{infcomp} expresses its compatibility with the polarization. In other words, the image of $\delta$ is the tangent space to the image, via the period map, of a variation of Hodge structure.

Given a variation of Hodge structure $(S, \bV, \cF^\bullet, \nabla)$ and a point $s\in S$, there is an associated infinitesimal variation of Hodge structure with $V_\bZ = \bV_s$, $T= T_s S$ and the map $\delta$ given by the differential of the period map of $(S, \bV, \cF^\bullet, \nabla)$ at $s$. A first interesting application of the theory of infinitesimal variations of Hodge structures was to curves: Griffiths observed \cite[Section 5]{Griffiths1983-Inf3} that if a curve $C$ is neither hyperelliptic, nor trigonal, nor a smooth plane quintic, then it is determined by its universal infinitesimal variation of Hodge structure. This follows from the fact that the kernel of the codifferential of the period map can be identified with the space of quadrics containing the canonical curve. Note that general curves of genus $\geq 5$ satisfy the above hypothesis. For curves of genus $3$, the period map is surjective hence the codifferential of the period map is $0$ and cannot determine the curve. We shall see below that in the cases of genus $3$ and $4$ one can use infinitesimal invariants of normal functions to determine the curve.

\section{The Griffiths infinitesimal invariant of a normal function}\label{SectInf}

Let $(S, \bV, \cF^\bullet, \nabla)$ be an integral variation of Hodge structure of odd weight $k= 2l-1$ on the complex manifold $S$, and let $\nu$ be a normal function on its associated family of intermediate jacobians. Griffiths defined an infinitesimal invariant associated to $\nu$ whose vanishing is a necessary condition for $\nu$ to be locally constant. Here we first define this invariant and its refinement by Green \cite{Green1989-Inf}, then we discuss some of its applications.

By Griffiths transversality and flatness, for all $p,r$, the connection $\nabla$ induces maps
\[
\onabla : \cV^{l-p,l+p-1} \otimes \Omega^r_S \lra \cV^{l-p-1, l+p} \otimes \Omega^{r+1}_S
\]
such that the sequence
\[
0 \lra \cV^{l-p, l+p-1} \lra \cV^{l-p-1, l+p} \otimes \Omega^1_S \lra \cV^{l-p-2, l+p+1} \otimes \Omega^2_S \lra \ldots
\]
is a complex. For any local lifting $\tnu \in \cV$ of $\nu$, we have
\[
\nabla (\tnu) \in \cF^{l-1} \otimes \Omega^1_S.
\]
One can show (see, e.g., \cite[Vol. II, Chapter 7]{Voisin2003-Hodge}) that the projection of $\nabla (\tnu)$ to $\cV^{l-1,l} \otimes \Omega^1_S$ lies in the kernel of the differential $\onabla$ and its class modulo the image $\Image (\onabla : \cV^{l,l-1} \ra \cV^{l-1,l} \otimes \Omega^1_S)$ depends only on $\nu$ and not on the choice of the lifting $\tnu$. We have the following

\begin{definition}
The infinitesimal invariant $\delta_1 (\nu)$ is the class of $\onabla (\tnu)$ in the cohomology group
\[
H^1 (\cV^{.,.}\otimes \Omega^._S) = \frac{\Ker (\onabla : \cV^{l-1,l} \otimes \Omega^1_S \ra \cV^{l-2,l+1} \otimes \Omega^2_S)}{\Image (\onabla : \cV^{l,l-1} \ra \cV^{l-1,l} \otimes \Omega^1_S)}.
\]
For $s\in S$, the infinitesimal invariant $\delta_1 (\nu)_s$ of $\nu$ at $s$ is the class of $\onabla (\tnu)_s$ in
\[
H^1 (\cV^{.,.}_s\otimes \Omega^._{S,s}) = \frac{\Ker (\onabla_s : \cV^{l-1,l}_s \otimes \Omega^1_{S,s} \ra \cV^{l-2,l+1}_s \otimes \Omega^2_{S,s})}{\Image (\onabla_s : \cV^{l,l-1}_s \ra \cV^{l-1,l}_s \otimes \Omega^1_{S,s})}.
\]
\end{definition}

Similarly, one has the complexes
\begin{equation}\label{eqfiltrationcomplex}
0 \lra \cF^l \lra \cF^{l-1} \otimes \Omega^1_S \lra \cF^{l-2} \otimes \Omega^2_S \lra \ldots.
\end{equation}
Mark Green \cite{Green1989-Inf} (and also Robert Bryant, see \cite[Appendix to Section 6(a)]{Griffiths1983-Inf3}) observed that one could define a slightly more general invariant whose vanishing is equivalent to $\nu$ being locally constant.
\begin{definition}
The infinitesimal invariant $\delta (\nu)$ is the class of $\onabla (\tnu)$ in the cohomology group
\[
H^1 (\cF^.\otimes \Omega^._S) = \frac{\Ker (\onabla : \cF^{l-1} \otimes \Omega^1_S \ra \cF^{l-2} \otimes \Omega^2_S)}{\Image (\onabla : \cF^{l} \ra \cF^{l-1} \otimes \Omega^1_S)}.
\]
For $s\in S$, the infinitesimal invariant $\delta (\nu)_s$ of $\nu$ at $s$ is the class of $\onabla (\tnu)_s$ in
\[
H^1 (\cF^._s\otimes \Omega^._{S,s}) = \frac{\Ker (\onabla_s : \cF^{l-1}_s \otimes \Omega^1_{S,s} \ra \cF^{l-2}_s \otimes \Omega^2_{S,s})}{\Image (\onabla_s : \cF^{l}_s \ra \cF^{l-1}_s \otimes \Omega^1_{S,s})}.
\]
\end{definition}

The Hodge filtration also induces filtrations on the complexes \eqref{eqfiltrationcomplex} as well as on their cohomology sheaves. The image of $\delta (\nu)$ in the first graded piece of this filtration is the invariant $\delta_1(\nu)$. Green also observed that if $\delta_1(
\nu)$ is zero, then the image $\delta_2(\nu)$ of $\delta (\nu)$ in the second graded piece of the filtration is well-defined. One could continue this to define infinitesimal invariants $\delta_1(\nu), \ldots, \delta_{l+1} (\nu)$ such that $\delta_p(\nu)$ is well-defined when $\delta_1(\nu) = \ldots = \delta_{p-1} (\nu) = 0$. Then the vanishing of $\delta (\nu)$ is equivalent to the vanishing of $\delta_i (\nu)$ for $i=1, \ldots, l+1$.

The non-vanishing of any of the infinitesimal invariants implies that the normal function and all of its multiples are not locally constant. In particular, the normal function is not torsion.

Infinitesimal invariants of normal functions have had many interesting applications to the study of cycles on algebraic varieties.

A first application by Griffiths was to show that a general curve $C$ of genus $4$ can be reconstructed from the infinitesimal invariant of the normal function with value the difference of the two $g^1_3$s on $C$ \cite[Section 6(d)]{Griffiths1983-Inf3}. More precisely, the infinitesimal invariant can be identified with a cubic polynomial whose trace on the quadric determined by the universal infinitesimal deformation of $C$ (which is the unique quadric containing the curve) is the canonical image of $C$.

As we mentioned in Section \ref{sectIntJ}, infinitesimal invariants are the main ingredient of the proof of the Green-Voisin theorem \cite{Green1989-Inf}.

Claire Voisin \cite{Voisin1992-GG} used infinitesimal invariants to give a second proof of Clemens' theorem that the Griffiths group of a generic quintic threefold is not finitely generated over $\bQ$.

The infinitesimal invariant was also the main tool used by Collino--Pirola \cite{CollinoPirola1995} and Collino-Naranjo-Pirola \cite{CollinoNaranjoPirola} to prove the results that we described in the introduction.

\section{Griffiths' formula for the infinitesimal invariant}\label{secGriffform}

In their computation of the infinitesimal invariant, Collino and Pirola \cite{CollinoPirola1995} used a formula of Griffiths \cite[Section 6(e)]{Griffiths1983-Inf3} for the first infnitesimal invariant $\delta_1 (\nu)$ of the normal function $\nu$ of a family of algebraic cycles on a family of smooth projective varieties. In this section we present Griffiths' formula.

Let $\cX \ra S$ be a family of smooth projective varieties of dimension $n=2m+1$, choose a point $s_0 \in S$ and put $X := X_{s_0}$. Via a $C^\infty$ trivialization $\cX \cong X \times S$, the varying complex structures on the fibers of $\cX \ra S$ are described by a family of operators $\odel_s : A^0 (X) \ra A^1 (X)$ where $A^p (X)$ is the space of global $C^\infty$ $p$-forms on $X$. Then
\[
\theta_s := \odel_s - \odel,
\]
where $\odel := \odel_{s_0}$ and $\{\theta_s\}$ is a holomorphic family of vector valued $(0,1)$-forms on $X$ which can be locally written as
\[
\theta_s(z) = \sum_{i,j} \theta^i_j (s,z) \frac{\del}{\del z_i} \otimes d \oz_j
\]
in terms of local holomorphic coordinates on $X$. Using the integrability condition
\[
\odel \theta_s - \frac{1}{2} [\theta_s , \theta_s ] =0,
\]
one verifies that, for any $\xi \in T_{s_0} S$, the derivative $\frac{\del \theta_s}{\del \xi}$ of $\theta_s$ in the direction of $\xi$ represents the image of $\xi$ under the Kodaira--Spencer map $T_{s_0} S \ra H^1 (X, T_X)$.

We have the variation of Hodge structure $(S, \bH^{2m+1} := R^{2m+1}\phi_* \bZ_\cM, \cF^\bullet, \nabla)$ and its associated family of intermediate jacobians
\[
\cJ^m_S := \frac{\cH^{2m+1}}{\cF^{m+1} \oplus \bH^{2m+1}}.
\]
Consider now a holomorphic family $Z_s \subset X_s$ of homologically trivial algebraic cycles of codimension $m+1$. We have the normal function
\[
\begin{array}{cccc}
\nu : & S & \lra & \cJ^m_S \\
 & s & \longmapsto & AJ_s (Z_s),
\end{array}
\]
where $AJ_s$ is the Abel-Jacobi map for $X_s$ (see Section \ref{sectIntJ}).

For any $s\in S$, cup-product, followed by integration gives a non-degenerate pairing
\[
\Psi_s : H^{m,m+1}(X_s) \otimes H^{m+1,m}(X_s) \lra \bC,
\]
which is also the restriction of the polarization on $H^{2m+1} (X_s)$. Pairing $\delta_1 (\nu)$ with a tangent vector $\xi \in T_{s_0} S$ gives $\nabla_\xi (\tnu) \in H^{m,m+1} (X)$ which is determined by the values $\Psi (\nabla_\xi (\tnu), \omega)$ for all $\omega \in H^{m+1,m} (X)$. Griffiths computed these values. Let $\eta_Z$ be a normal vector field on $Z$ which projects injectively to $S$ and, for $\omega \in H^{m+1,m}(X)$, denote $i (\eta_Z)\omega$ the contraction of $\omega$ against $\eta_Z$. Let $\zeta$ be a differential form of type $(m,m)$ on $X$ such that $\odel (\zeta) = \frac{\del \theta_s}{\del \xi} \cup \omega$.

\begin{lemma}(Griffiths \cite[Section 6]{Griffiths1983-Inf3})
With the above notation, we have
\[
\Psi (\nabla_\xi (\tnu), \omega) = \int_Z \zeta - \int_Z i(\eta_Z) \omega
\]
\end{lemma}

\section{Voisin's formula for the infinitesimal invariant}\label{sectvoisin}

For their generalization of the results of Collino--Pirola \cite{CollinoPirola1995}, Pirola and Zucconi \cite{PirolaZucconi2003} used a formula of Voisin \cite{Voisin1988-Inf} for the computation of the first infinitesimal invariant. In this section we present Voisin's formula.

Let $\tX$ be a smooth projective variety of dimension $n+1 = 2m+2$ and let $X \subset \tX$ be a smooth divisor, moving in a linear system $L$. Let $Z$ be an $(m+1)$-dimensional algebraic cycle on $\tX$ whose restriction to $X$ is homologically trivial. In other words, the cohomology class $[Z]$ is primitive, i.e.,
\[
[Z]\in H^{m+1} (\tX, \Omega_{\tX}^{m+1} )^{prim}
\]
where the superscript $prim$ indicates the kernel of the restriction map:
\[
H^{m+1} (\tX, \Omega_{\tX}^{m+1} )^{prim} := \Ker \left( H^{m+1} (\tX, \Omega^{m+1}_{\tX} ) \stackrel{res}{\lra} H^{m+1} (X, \Omega^{m+1}_{X} ) \right).
\]
Then, on the locus $L_{reg}$ in $L$ parametrizing smooth divisors, we have the family of divisors $\cX \ra L_{reg}$, the variation of Hodge structure $(L_{reg}, \bH^{2m+1} := R^{2m+1}\phi_* \bZ_\cX, \cF^\bullet, \nabla)$ and its associated family of intermediate jacobians
\[
\cJ^m_{L_{reg}} := \frac{\cH^{2m+1}}{\cF^{m+1} \oplus \bH^{2m+1}}.
\]
The cycle $Z$ restricts to a homologically trivial cycle $Z_s$ on every fiber $X_s$ of $\cX \ra L_{reg}$ and we have the associated normal function
\[
\begin{array}{cccc}
\nu : & {L_{reg}} & \lra & \cJ^m_{L_{reg}} \\
 & s & \longmapsto & AJ_s (Z_s).
\end{array}
\]
For $s\in L_{reg}$, the first infinitesimal invariant $\delta_1 (\nu)$ is an element of
\[
\begin{split}
H^1 (\cH^{.,.}_s\otimes \Omega^._{L_{reg},s}) = \frac{\Ker (\onabla_s : H^{m+1} (X_s, \Omega^m_{X_s}) \otimes \Omega^1_{L_{reg},s} \ra H^{m+2} (X_s, \Omega^{m-1}_{X_s}) \otimes \Omega^2_{L_{reg},s})}{\Image (\onabla_s : H^m (X_s, \Omega^{m+1}_{X_s}) \ra H^{m+1} (X_s, \Omega^m_{X_s}) \otimes \Omega^1_{L_{reg},s})} \\
\subset \Coker \left(\onabla_s : H^m (X_s, \Omega^{m+1}_{X_s}) \ra H^{m+1} (X_s, \Omega^m_{X_s}) \otimes \Omega^1_{L_{reg},s}\right).
\end{split}
\]
Dually, the infinitesimal invariant is an element of
\[
\Ker \left( \prescript{t}{}{\onabla_s} :H^m (X_s, \Omega^{m+1}_{X_s}) \otimes T_{L_{reg},s} \ra H^{m+1} (X_s, \Omega^m_{X_s})\right)^\vee,
\]
where, with the notation of the previous section, we are using the pairing $\Psi$ to identify $H^{m+1, m} (X_s) = H^m (X_s, \Omega^{m+1}_{X_s})$ with the dual of $H^{m,m+1} (X_s) = H^{m+1} (X_s, \Omega^m_{X_s})$. The map $\onabla$ and its transpose $\prescript{t}{}{\onabla}$ are both given by cup-product via the Kodaira--Spencer map. This means, for instance, that for $\omega \in H^{m+1, m} (X_s)$ and $v\in T_{L_{reg},s}$, we have
\[
\prescript{t}{}{\onabla} (\omega \otimes v) = \omega \cup \xi_v
\]
where $\xi_v \in H^1 (X_s, T_{X_s})$ is the image of $v$ via the Kodaira--Spencer map.

Identifying the tangent space to $L$ at $s$ with the quotient $H^0 (X, L) / \langle f_s \rangle$, where we denote $f_s$ an element of $H^0 (X, L)$ with divisor of zeros $X_s$, the Kodaira--Spencer map for $L_{reg}$ is induced by the connecting homomorphism of the usual normal bundle sequence
\[
0 \lra T_{X_s} \lra T_X |_{X_s} \lra \cO_{X_s} (L) \lra 0,
\]
using the natural embedding $H^0 (X, L) / \langle f_s \rangle \subset H^0 (X_s, \cO_{X_s} (L))$.

Put
\[
H^{m+1} (\tX, \Omega_{\tX}^{m+1} |_{X_s})^{prim} := \Ker \left( H^{m+1} (\tX, \Omega^{m+1}_{\tX} |_{X_s}) \stackrel{res}{\lra} H^{m+1} (X_s, \Omega^{m+1}_{X_s} )\right).
\]
It follows that we have the commutative diagram with exact rows
\[
\xymatrix{
& & H^{m+1} (\tX, \Omega_{\tX}^{m+1} )^{prim} \ar[d]^{res} \\
H^m (X_s, \Omega^{m+1}_{X_s}) \ar@{=}[d] \ar[r] & H^{m+1} (X_s, \Omega^m_{X_s}(-L)) \ar[d] \ar[r] & H^{m+1} (\tX, \Omega_{\tX}^{m+1} |_{X_s})^{prim} \ar[d]^{\oalpha} \\
H^m (X_s, \Omega^{m+1}_{X_s}) \ar[r]^{\onabla_s\phantom{mnmnmnmm}} & H^{m+1} (X_s, \Omega^m_{X_s}) \otimes (H^0 (X, L) / \langle f_s \rangle)^\vee \ar[r] & \Coker (\onabla_s),}
\]
where the middle row is part of the cohomology sequence of the sequence
\[
\xymatrix{
0 \ar[r] & \Omega^m_{X_s}(-L) \ar[r] & \Omega^{m+1}_{\tX} |_{X_s} \ar[r] & \Omega^{m+1}_{X_s} \ar[r] & 0.}
\]
Voisin proved the following \cite{Voisin1988-Inf}.

\begin{proposition}
With the map $\oalpha$ defined by the above diagram, we have
\[
\delta_1(\nu_Z) = \oalpha (res ([Z])).
\]
\end{proposition}
When $L$ is a base point free pencil, we have $\cO_{X_s} (L) \cong \cO_{X_s}$ and hence the exact sequence
\[
\xymatrix{
0 \ar[r] & \Omega^{m+1}_{\tX} \ar[r] & \Omega^{m+1}_{\tX} (L) \ar[r] & \Omega^{m+1}_{\tX} |_{X_s} \ar[r] & 0,}
\]
whose connecting homomorphism
\[
\xymatrix{
\prescript{t}{}{res} : H^m (X_s, \Omega^{m+1}_{\tX} |_{X_s} ) \ar[r] & H^{m+1} (\tX, \Omega^{m+1}_{\tX})}
\]
is the transpose of $res$.
We also have the exact sequence
\[
\xymatrix{
0 \ar[r] & \Omega^m_{X_s} \ar[r] & \Omega^{m+1}_{\tX} |_{X_s} \ar[r] & \Omega^{m+1}_{X_s} \ar[r] & 0,}
\]
whose connecting homomorphism can be identified with $\prescript{t}{}{\onabla}$ after a choice of basis for $T_s L \cong \bC$.

Given $\omega \otimes v \in \Ker \prescript{t}{}{\onabla}$, let $\gamma$ be a pre-image for $\omega \otimes v$ in $H^m (X_s, \Omega^{m+1}_{\tX} |_{X_s})$. Voisin also proves that, dually, we have
\[
\delta_1 (\nu_Z) (\omega\otimes v) = \prescript{t}{}{res} (\gamma) \cup [Z].
\]

\section{The Collino--Pirola--Zucconi adjunction map}\label{defadj}

Collino--Pirola \cite{CollinoPirola1995} introduced an adjunction map for infinitesimal deformations of smooth varieties that they used, in conjunction with Griffiths' formula in Section \ref{secGriffform}, to compute certain values of the infinitesimal invariant of some normal functions. We present here a combination of results contained in \cite{CollinoPirola1995}, \cite{CollinoNaranjoPirola}, \cite{PirolaZucconi2003}.

Suppose given a smooth projective $n$-dimensional variety $X$ with a locally free sheaf $\cF$ of rank $n$ and an extension
\[
0 \lra \cO_X \lra \cE \lra \cF \lra 0
\]
with extension class $\xi\in \Ext^1 (\cF , \cO_X) \cong H^1 (X, \cF^*)$. Denote $\del_\xi : H^0 (X, \cF) \ra H^1(X, \cO_X)$ the first connecting homomorphism of the long exact sequence of cohomology associated to the sequence above. Put $m := h^0 (X, \cF)$, assume $m >n$, and let $\bG = G(n+1, m)$ be the Grassmannian of $n+1$ dimensional subspaces of $H^0 (X, \cF)$ with universal subbundle $\cU$. Let
\[
\Gamma := \{ (W, \xi ) \mid W \subset \Ker \del_\xi \} \subset \bG \times \bP \Ext^1 (\cF, \cO_X)
\]
be the incidence correspondence. Let $\cU_\Gamma$ be the pull-back of $\cU$ to $\Gamma$ via the first projection and let $\cO_\Gamma (1)$ be the pull-back of $\cO_{\bP \Ext^1 (\cF, \cO_X)} (1)$ via the second projection. The natural map $\Lambda^n H^0 (X, \cF) \ra H^0 (X, \Lambda^n \cF)$ gives rise to the natural map
\[
\Lambda^n \cU_\Gamma \lra H^0 (X, \Lambda^n \cF) \otimes \cO_\Gamma
\]
whose cokernel we denote $\cG$. The adjunction map is the natural morphism
\[
\alpha : \Lambda^{n+1}\cU_\Gamma \otimes \cO_\Gamma (-1) \lra \cG
\]
which can be concretely described as follows. Given a subspace $W \subset \Ker \del_\xi$ of dimension $n+1$, choose a basis $s_1, \ldots, s_{n+1}$ of $W$ which gives an element $s_1 \wedge \ldots \wedge s_{n+1}$ of $\Lambda^{n+1} W$. Lift the elements $s_i$ to elements $s_i'$ of $H^0 (X, \cE)$. The class $\alpha (s_1 \wedge \ldots \wedge s_{n+1})$ is the image of $s_1' \wedge \ldots \wedge s_{n+1}'$ via the composition
\[
\Lambda^{n+1} H^0 (X, \cE) \lra H^0 (X, \Lambda^{n+1} \cE) \cong H^0 (X, \Lambda^n \cF) \surj H^0 (X, \Lambda^n \cF) / \Lambda^n W.
\]
It is not difficult to check that this is indeed a well-defined map. It induces the above bundle map $\alpha$.

\section{The infinitesimal invariant and the adjunction map}\label{sectlink}

Given a commutative diagram
\[
\xymatrix{
\cX \ar[dr]_{\pi} \ar[rr]^{\Phi} & & \cA \ar[dl]^{\rho} \\
& B &}
\]
where $\cX$ is a family of smooth projective varieties of dimension $n$, $\cA$ is a family of abelian varieties of dimension $a$ over the smooth analytic base $B$, we have the natural cycle
\[
\cZ := [\cY] - [\cY]^-
\]
where $\cY$ is the image of $\cX$ via $\Phi$, $[\cY]$ is the associated cycle on $\cA$ and $[\cY]^-$ is its image under $-\Id$. For $b\in B$, we denote the fibers of $\cX, \cY, \cZ, \cA$, by $X_b, Y_b, Z_b, A_b$ respectively.

After possibly shrinking $B$, we can replace $\cX$ by a desinguarization of $\cY$ and replace $\cA$ by the family of its abelian subvarieties generated by the fibers of $\cY$ so that the following conditions are satisfied.
\begin{enumerate}
\item For all $b\in B$, the restriction $\phi_b : X_b \ra A_b$ is of degree $1$ onto its image $Y_b$, and
\item for all $b\in B$, the image $Y_b$ generates $A_b$ as a group.
\end{enumerate}
Families with the two properties above were called ``of Albanese type'' in \cite{PirolaZucconi2003}.

The cycle $\cZ$ is homologically trivial of codimension $a-n$ in each fiber. It therefore has an associated normal function $\nu_\cZ$:
\[
\begin{array}{rcl}
\nu_\cZ : B & \lra & \cJ^{a-n} \\
b & \longmapsto & AJ (Z_b - Z_b^-)
\end{array}
\]
where
\[
\cJ^{a-n}_B :=  \frac{\cH^{2a-2n-1}(\cA)}{\cF^{a-n} \oplus \bH^{2a-2n-1}}
\]
is the family of $(a-n)$-th intermediate jacobians of $\cA$.

Suppose now that $B$ is a disc with center $b$ and let $v\in T_b B$ be a nonzero vector tangent to $B$. The usual cotangent bundle sequence
\[
\xymatrix{
0 \ar[r] & \cO_{A_b} \ar[r] & \Omega^1_{\cX} |_{X_b} \ar[r] & \Omega^1_{X_b} \ar[r] & 0}
\]
has extension class $\xi_v \in H^1 (X_b, T_{X_b}) = \Ext^1 (\Omega^1_{X_b}, \cO_{A_b})$, the Kodaira-Spencer image of $v$. We use the adjunction map for this extension of locally free sheaves.

Assume that there exist $n+1$ global sections $s_1, \ldots , s_{n+1}$ of $\Omega^1_\cA$ whose pull-backs $\eta_1, \ldots , \eta_{n+1}$ to $X$ lie in the kernel of the connecting homomorphism $\del_{\xi_v}$. Put $\Omega := s_1 \wedge \ldots \wedge s_{n+1}$ and let $\omega$ be a form representing the class $\alpha (\eta_1 \wedge \ldots \wedge \eta_{n+1})\in H^0 (X, \Omega^n_X) / \Lambda^n W$, where $W$ is the span of $\eta_1, \ldots , \eta_{n+1}$ in $H^0 (X, \Omega^1_X)$. Let
\[
\xymatrix{
\gamma : H^n (A_b, \Omega^{n+1}_\cA |_{A_b} ) \ar[r] & H^n (A_b, \Omega^{n+1}_{A_b})}
\]
be the restriction map.

Pirola and Zucconi \cite{PirolaZucconi2003} prove
\begin{proposition}
For all $\sigma \in H^0 (A_b, \Omega^n_{A_b})$, we have
\begin{enumerate}
\item $\gamma (\Omega \otimes \osigma) \otimes v \in \Ker \prescript{t}{}{\onabla_b}$,
\item $\delta_1 (\nu_\cZ) (\gamma (\Omega \otimes \osigma) \otimes v) = 2 \int_X \omega \wedge \overline{\phi_b^* \sigma}$.
\end{enumerate}
\end{proposition}

This formula was one of the main tools used by Pirola and Zucconi in \cite{PirolaZucconi2003} and by Collino, Naranjo and Pirola \cite{CollinoNaranjoPirola} to prove the results that we described in the introduction. It was originally proved by Collino and Pirola \cite{CollinoPirola1995} in the case where $\cX \ra \cA$ is a family of curves in their jacobians, using the formula of Griffiths in Section \ref{secGriffform}.

In the case of curves in their jacobians, the relevant intermediate jacobian is $J^{g-1}(JX) = \frac{H^{2g-3} (JX, \bC)}{F^{g-1} + H^{2g-3}(JX, \bZ)} = \frac{F^2H^3(JX, \bC)^*}{H_3(JX, \bZ)}$. This intermediate jacobian always contains a nontrivial abelian subvariety which is generated by translating the Abel-Jacobi image of a copy of $X$ in $JX$. One must therefore first quotient by this abelian subvariety: this quotient is the primitive intermediate jacobian $P^{g-3} (JX)$. To show that the Ceresa cycle is not algebraically equivalent to zero, Collino and Pirola show that the normal function of the Ceresa cycle is generically non-abelian, i.e., it does not take values in the abelian part of $P^{g-3} (JX)$. In fact, using the above calculation, they show that the codimension of any subspace $B$ of $\cM_g$ where the normal function of $X - X^-$ is abelian is greater than $\frac{g+2}{3}$ if $g\geq 4$ and greater than $2$ if $g=3$ and $B$ is not contained in the hyperelliptic locus.
Complementing Griffiths' results (see Sections \ref{secInfVar} and \ref{SectInf}), they also show that the infinitesimal invariant can be identified with the equation of the canonical curve when the genus is $3$.







\providecommand{\bysame}{\leavevmode\hbox to3em{\hrulefill}\thinspace}
\providecommand{\MR}{\relax\ifhmode\unskip\space\fi MR }
\providecommand{\MRhref}[2]{%
  \href{http://www.ams.org/mathscinet-getitem?mr=#1}{#2}
}
\providecommand{\href}[2]{#2}

\end{document}